\documentclass{amsart}
\usepackage{amssymb}
\usepackage{amsfonts,amsmath,amssymb,mathrsfs,verbatim}
\usepackage[dvips]{graphicx}
\usepackage{graphics}
\usepackage{psfrag}
\usepackage{tcolorbox}

\definecolor{UBCblue}{rgb}{0.24706, 0.63725, 0.46667} 

\begin{document}
\newtheorem{thm1}{Theorem}[section]
\newtheorem{lem1}[thm1]{Lemma}
\newtheorem{rem1}[thm1]{Remark}
\newtheorem{def1}[thm1]{Definition}
\newtheorem{cor1}[thm1]{Corollary}
\newtheorem{defn1}[thm1]{Definition}
\newtheorem{prop1}[thm1]{Proposition}
\newtheorem{ex1}[thm1]{Example}
\newtheorem{alg1}[thm1]{Algorithm}


\title[toric ideals of graphs]{An algorithm for computing the universal Gr\"obner Basis of graph ideals}
\author{Yannis C. Stamatiou}
\address{Dept. of Business Administration, University of Patras
and Computer Technology Institute and Press - ``Diophantus'', Patra,
Greece}
\email{stamatiu@ceid.upatras.gr}
\author{Christos Tatakis}
\address{Department of Mathematics, University of Ioannina,
Ioannina 45110, Greece}
\email{chtataki@cc.uoi.gr}
\thanks{}

\subjclass[2000]{Primary  14M25, 68Q40, 68W30, 13P10, 05C85}

\date{}

\dedicatory{}

\begin{abstract}
\par
The universal Gr\"obner basis of an ideal is a Gr\"obner basis with respect to all term orders simultaneously. The aim of this paper is to present an algorithmic approach to compute the universal Gr\"obner basis for the toric ideal corresponding
to an undirected graph, based on the theoretically knowledge of this set and on
 a recent, efficiently computable
algorithmic characterization of the Graver basis of the ideal.
\end{abstract}

\maketitle

\section{Introduction}
\par
A Gr\"{o}bner basis is a specific generating set of an ideal $I$
over a polynomial ring $\mathbb{K}[x_1,\ldots,x_n]$, where $\mathbb{K}$ is a field. It has extremely useful algebraic properties
and it is relatively easy to extract information about the ideal,
given its Gr\"{o}bner basis. The study of the Gr\"{o}bner
bases has become a major research topic in commutative algebra,
combinatorics and computer science. Gr\"obner basis theory provides
the foundations for many algorithms in commutative algebra and
algebraic geometry. Its importance stems from its wide applicability
in problems coming from diverse disciplines such as mathematics,
combinatorics, computer design theory, symbolic computation,
 integer programming,
engineering, computer technology and cryptography. The concept of a Gr\"{o}bner basis was introduced by Buchberger at
1965 who named this set after his supervisor W. Gr\"obner.
Buchberger's algorithm is the most well known algorithmic method for
computing a Gr\"obner basis for an ideal $I$ of a polynomial ring.
Most of the symbolic computation software packages, such as CoCoA,
Macaulay, Mathematica, Maple and Singular, include algorithms for
computing this set. For more information on Gr\"obner bases and
corresponding algorithms, see \cite{AL,BW,CLS}. \par
The universal Gr\"{o}bner basis is defined as the union of all reduced Gr\"{o}bner
bases $G_<$ of the ideal $I$, as $<$ runs over all term orders. It is a finite subset of $I$ and it is a
Gr\"{o}bner basis for the ideal with respect to all admissible term
orders (for more see \cite{St}). They were introduced by
V.~Weispfenning \cite{W} and N.~Schwartz \cite{S}. There are not
many classes of ideals for which we know their universal Gr\"{o}bner
bases. In general, characterizing and computing this set is a very difficult and computationally demanding
problem. One of the known classes, is the class of the toric ideals associated to incidence matrices of graphs. The
structure of the universal Gr\"{o}bner basis of a toric ideal of a
graph $G$, was characterized theoretically in \cite[Theorem 3.4.]{TT2}.\par  The goal of this paper is to
transform this theorem into an algorithm in order to compute this set. Computing the universal Gr\"obner basis of an ideal, is in general extremely complicated, despite the existence of a few algorithms, valuable from a theoretical point of view but of a limited practical use. One would question the necessity of computing the universal Gr\"obner basis of the ideal, while at the same time we could just use a Gr\"obner basis of it under fixed term order, which would also be proper and there are already available algorithms for it. The reply stems from the fact that in certain problems one needs to compute Gr\"obner bases for the same  
toric ideal but each time with different term order. For  
example, problems coming from integer programming where the term  
order depends on a cost function that should be minimized. Then,  
instead of a Gr\"obner computation each time that cost function changes,  
if the universal Gr\"obner basis is known, we just have to compare  
degrees in every binomial in the universal Gr\"obner basis to produce
a Gr\"obner basis for the toric ideal. \par
The structure of this paper is as follows. In Section 2 we present the terminology and the theory which someone needs to follow. The algorithm is presented in Section 3. Its input is the Graver basis of the ideal, which it is known algorithmically from \cite{OgTak}. We compute the complexity of the algorithm, see Theorem \ref{thm}, and we prove that the algorithm decides in polynomial time whether an element of the Graver basis of
the ideal belongs also to its universal Gr\"obner basis, see Corollary \ref{cor-compl}.

\section{Graver basis and the universal Gr\"obner basis of a graph ideal}
\par  In this section we define the toric ideals of graphs and next we give some basic elements of graph theory
which will be useful in the description of the Graver and universal
Gr\"obner basis of a toric ideal of a graph $G$. \par
 Let $A=\{\textbf{a}_1,\ldots,\textbf{a}_m\}\subseteq \mathbb{N}^n$
be a vector configuration in $\mathbb{Q}^n$ and
$\mathbb{N}A:=\{l_1\textbf{a}_1+\cdots+l_m\textbf{a}_m \ | \ l_i \in
\mathbb{N}\}$ the corresponding affine semigroup.  We grade the
polynomial ring $\mathbb{K}[x_1,\ldots,x_m]$ over an arbitrary field
$\mathbb{K}$ by the semigroup $\mathbb{N}A$ setting
$\deg_{A}(x_i)=\textbf{a}_i$ for $i=1,\ldots,m$. For
$\textbf{u}=(u_1,\ldots,u_m) \in \mathbb{N}^m$, we define the
$A$-\emph{degree} of the monomial
$\textbf{x}^{\textbf{u}}:=x_1^{u_1} \cdots x_m^{u_m}$ to be \[
\deg_{A}(\textbf{x}^{\textbf{u}}):=u_1\textbf{a}_1+\cdots+u_m\textbf{a}_m
\in \mathbb{N}A.\]  The \emph{ toric ideal} $I_{A}$ associated to
$A$ is the prime ideal generated by all the binomials
$\textbf{x}^{\textbf{u}}- \textbf{x}^{\textbf{v}}$ such that
$\deg_{A}(\textbf{x}^{\textbf{u}})=\deg_{A}(\textbf{x}^{\textbf{v}})$,
see \cite{St}. For such binomials, we set
$\deg_A(\textbf{x}^{\textbf{u}}-
\textbf{x}^{\textbf{v}}):=\deg_{A}(\textbf{x}^{\textbf{u}})$. An
irreducible binomial $\textbf{x}^{\textbf{u}}-
\textbf{x}^{\textbf{v}}$ in $I_A$ is called \emph{ primitive} if
there exists no other binomial
 $\textbf{x}^{\textbf{w}}- \textbf{x}^{\textbf{z}}$ in $I_A$
such that $\textbf{x}^{\textbf{w}}$ divides $
\textbf{x}^{\textbf{u}}$ and $\textbf{x}^{\textbf{z}}$ divides $
\textbf{x}^{\textbf{v}}$. The set of primitive binomials forms the
Graver basis of $I_A$ and is denoted by $Gr_A$.\par
We study a special class of toric ideals which arise from graphs. Let $G$ be a connected, undirected, finite, simple graph on the vertex set
$V(G)=\{v_{1},\ldots,v_{n}\}$.  Let $E(G)=\{e_{1},\ldots,e_{m}\}$ be
the set of edges of $G$ and $\mathbb{K}[e_{1},\ldots,e_{m}]$
 the polynomial ring in the $m$ variables $e_{1},\ldots,e_{m}$ over an arbitrary field $\mathbb{K}$.  We
will associate each edge $e=\{v_{i},v_{j}\}\in E(G)$  with the element
$a_{e}=v_{i}+v_{j}$ in the free abelian group $ \mathbb{Z}^n $, where $v_{i}=(0,\ldots,0,1,0,\ldots,0)$ is the vector with 1 in the $i-$th coordinate. By $I_{G}$ we denote
 the toric ideal $I_{A_{G}}$ in
$\mathbb{K}[e_{1},\ldots,e_{m}]$, where  $A_{G}=\{a_{e}\ | \ e\in E(G)\}\subset \mathbb{Z}^n $.\par A \emph{walk} connecting
$v_{i_{1}}\in V(G)$ and $v_{i_{s+1}}\in V(G)$ is a finite sequence
of the form
$$w=(\{v_{i_1},v_{i_2}\},\{v_{i_2},v_{i_3}\},\ldots,\{v_{i_s},v_{i_{s+1}}\})$$
with each $e_{i_j}=\{v_{i_j},v_{i_{j+1}}\}\in E(G)$. The
\emph{length} of the walk $w$ is the number $s$ of edges of the
walk. An even (respectively odd) walk is a walk of \emph{even}
(respectively odd) length. A walk
$w=(\{v_{i_1},v_{i_2}\},\{v_{i_2},v_{i_3}\},\ldots,\{v_{i_s},v_{i_{s+1}}\})$
is called \emph{closed} if $v_{i_{s+1}}=v_{i_1}$. A \emph{cycle} is
a closed walk
$$(\{v_{i_1},v_{i_2}\},\{v_{i_2},v_{i_3}\},\ldots,\{v_{i_s},v_{i_{1}}\})$$
with $v_{i_k}\neq v_{i_j},$ for every $ 1\leq k < j \leq s$. Note
that, although the graph $G$ has no multiple edges, the same edge
$e$ may appear more than once in a walk. In this case $e$ is called
a \emph{multiple edge of the walk $w$}. A graph is called biconnected if it is connected and does not contain a vertex whose removal increases the number of connected components of the graph. Also, we define the degree of a vertex of a graph (cor. of a walk) as the usual definition, which means as the number of the edges of the graph (cor. of the walk) incident to the vertex. Given an even closed walk
of the graph $G$, $w =(e_{i_1}, e_{i_2},\dots, e_{i_{2q}})$ we
denote by $B_w$ the binomial $$B_w=\prod
_{k=1}^{q} e_{i_{2k-1}}-\prod _{k=1}^{q}e_{i_{2k}}.$$ We remark that
$B_w\in I_G$. Moreover, it is known that the toric ideal $I_G$ is
generated by binomials of this form, see \cite{Vi}. We note that the same walk can be written in different ways but the corresponding binomials may differ only in the sign. Also note that for certain even closed walks $w$ the binomial $B_w$ may be zero, for example take $\zeta$ to be any walk then for the even
closed walk $w=(\zeta,-\zeta)$ we have $B_w = 0$.
\par

The relation between the Graver basis and the universal Gr\"{o}bner
basis, which is denoted by $\mathcal{ U}_A$, for a toric ideal $I_A$
was described by B.~Sturmfels:

\begin{prop1}\cite{St} For any toric ideal $I_A$ we have $\mathcal{ U}_A \subset  Gr_A$.
\end{prop1}

Therefore the knowledge of the Graver basis for
a toric ideal plays a key role for computing the
$\mathcal{ U}_A$. In order to describe the above sets in graph ideals, first, we give some extra graph notations. \par
We denote by $\textbf{w}$ the subgraph of $G$ with vertices the
vertices of the walk and edges the edges of the walk $w$. If $W$ is
a subset of the vertex set $V(G)$ of $G$ then the {\em induced
subgraph} of $G$ on $W$ is the subgraph of $G$ whose vertex set is
$W$ and whose edge set is $\{\{v, u\}\in E(G)|v,u\in W\}$. The walk
$w$ is primitive if and only if the binomial $B_w$ is primitive. A
\emph{cut edge} (respectively \emph{cut vertex}) is an edge
(respectively vertex) of the graph whose removal increases the
number of connected components of the remaining subgraph. A
\emph{block} is a maximal connected subgraph of a given graph $G$
which does not contain a cut vertex.\par The next corollary describes the elements $B_w\in I_G$ that belong to the Graver basis.

\begin{cor1}\cite[Corollary 3.3]{RTT} \label{primitive-graph}
Let $G$ be a graph and $W$ a connected subgraph of $G$. The subgraph
$W$ is the graph  ${\bf w}$ of a primitive walk $w$ if and only if
\begin{enumerate}
  \item  $W$ is an even cycle or
  \item  $W$ is not biconnected and
\begin{enumerate}
  \item every block of $W$ is a cycle or a cut edge and
  \item every cut vertex of $W$ belongs to exactly
  two blocks and separates the graph in two parts, the total number of
edges of the cyclic blocks in each part is odd.
\end{enumerate}
\end{enumerate}
\end{cor1}

In \cite{TT2}
the authors gave a necessary and sufficient
characterization of the universal Gr\"obner basis of $I_G$. In order to describe it, we need the notion of the pure cyclic block.
\par

Every even primitive walk $w=(e_{i_1},\ldots,e_{i_{2k}})$ partitions
the set of edges in the two sets $\textbf{w}^+= \{e_{i_j}|j \
\textit{ odd}\}, \textbf{w}^-=\{e_{i_j}|j \ \textit{ even}\}$,
otherwise the binomial $B_w$ would not be irreducible. The edges of
$\textbf{w}^+$ are called odd edges of the walk and those of
$\textbf{w}^-$ even. Note that for an even closed walk whether an
edge is even or odd depends on the edge that we start counting from.
Thus, it is not important to identify whether an edge is even or odd
but to separate the edges into two disjoint classes.

\begin{def1}\cite[Definition 3.1]{TT2} A cyclic
block $B$ of a primitive walk $w$ is called pure if all edges of $B$
are either in $\textbf{w}^+$ or in  $\textbf{w}^-$. A primitive walk
$w$ is called mixed if no cyclic block of $w$ is pure.
\end{def1}

The following Theorem describes the elements of the universal
Gr\"{o}bner basis of $I_G$, for any undirected graph $G$.

\begin{thm1}\cite[Theorem 3.4]{TT2} \label{UGB}
Let $w$ be a primitive walk. $B_w$ belongs to the universal
Gr\"{o}bner basis of $I_G$ if and only if $w$ is mixed.
\end{thm1}

Let see an example.

\begin{ex1}\label{exGraver}{\rm
We consider the graph $G$ depicted in Figure \ref{example} and let $$w=\{e_1,e_2,e_3,e_4,e_5,e_6,e_7,e_8,e_9,e_{10},e_{11},e_{12}\}$$
be an even walk which we can be seen in Figure \ref{example}. Let
$$B_w=\overbrace{e_1e_3e_5e_7e_9e_{11}}^{\textbf{w}^+}-\overbrace{e_2e_4e_6e_8e_{10}e_{12}}^{\textbf{w}^-}$$
be its corresponding binomial.

\begin{figure}[h]
\begin{center}
\psfrag{A}{$e_1$}\psfrag{B}{$e_2$}\psfrag{C}{$e_3$}\psfrag{D}{$e_4$}\psfrag{E}{$e_5$}\psfrag{F}{$e_6$}\psfrag{G}{$e_7$}\psfrag{H}{$e_8$}
\psfrag{I}{$e_9$}\psfrag{J}{$e_{10}$}\psfrag{K}{$e_{11}$}\psfrag{L}{$e_{12}$}
\psfrag{M}{$v_1$}\psfrag{N}{$v_2$}\psfrag{O}{$v_3$}\psfrag{P}{$v_4$}\psfrag{Q}{$v_5$}\psfrag{R}{$v_6$}\psfrag{S}{$v_7$}\psfrag{T}{$v_8$}\psfrag{U}{$v_9$}
\includegraphics{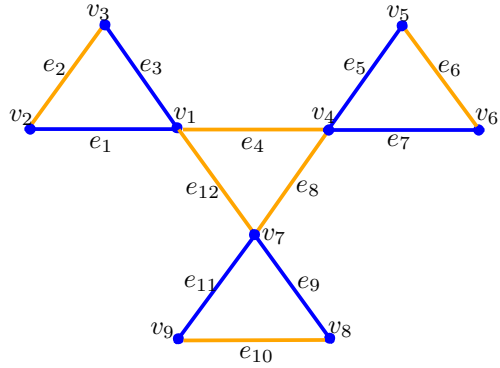}
\caption{The graph G}
\label{example}
\end{center}
\end{figure}

 The above graph $W$ is not biconnected, each of its blocks is a cycle and every cut vertex of $W$ belongs to exactly two blocks and
separates the graph in two parts, and the total number of edges of
the cyclic blocks in each part is odd. From Corollary \ref{primitive-graph} the walk $w$ is primitive and for the corresponding binomial $B_w$ it holds that $B_w\in Gr_G$.
Also $w$ is not mixed because of the existence of the block
$B_1=\{e_4,e_8,e_{12}\}$ which is not pure, since all of its edges
$\{e_4,e_8,e_{12}\}$ belong to $\textbf{w}^-$. Therefore from
Theorem~\ref{UGB}, we have that the binomial $B_w\notin \mathcal{ U}_G$. \par

Based on Corollary \ref{primitive-graph} one can check that the Graver basis of the toric ideal $I_G$ consists of the following ten elements:

$Gr_G=\{B_1=e_1e_3e_8-e_2e_4e_{12},B_2=e_4e_6e_8-e_5e_7e_{12},B_3=e_8e_{12}e_{10}-e_4e_9e_{11},B_4=e_1e_3e_5e_7-e_2e_4^2e_6,B_5=e_1e_3e_8^2e_6-e_2e_{12}^2e_7e_5,B_6=e_1e_3e_{11}e_9-e_2e_{12}^2e_{10},B_7=e_1e_3e_8^2e_{10}-e_2e_4^2e_9e_{11},B_8=e_{11}e_9e_4^2e_6-e_{10}e_{12}^2e_5e_7,B_9=e_{11}e_9e_7e_5-e_{10}e_8^2e_6,B_{10}=B_w=e_1e_3e_5e_7e_9e_{11}-e_2e_4e_6e_8e_{10}e_{12}\}$.

}
\end{ex1}

\section{Algorithmic description of the Universal Gr\"obner Basis of $I_G$}

In this section we state the main result of this paper which gives
an algorithmic description of the universal Gr\"obner basis 
of a toric ideal of a graph $G$. We remind that $A=\{\textbf{a}_1,\ldots,\textbf{a}_m\}\subseteq \mathbb{N}^n$
is a vector configuration in $\mathbb{Q}^n$. For the algorithm that we will present, its input is the set of the
primitive elements of a toric ideal of a graph $G$, i.e. its Graver basis. There are a lot of algorithms which give the Graver basis of an ideal. In  \cite{St}, Sturmfels gave algorithms for computing a Graver basis by Lawrence lifting of $A$ and the universal Gr\"obner basis for toric ideals with an algebraic geometry view. Moreover, he gave an algorithm, which takes as input  the universal Gr\"obner basis of $I_A$ in order to produce the State polytope of  the ideal, see \cite[Theorem 7.15]{St}.  In \cite{OgTak} M. Ogawa, H.
Hara and A. Takemura gave an algorithm for sampling elements from
the Graver basis set of $I_G$, which was associated with a simple
undirected graph for testing the beta model of graphs by Markov
chains, based on Monte Carlo methods. In general, the Graver basis
of toric ideals of graphs can be computed by symbolic computation software packages
such as the 4ti2 (see \cite{ti}) or Macaulay 2 (see \cite{GStil}).

The next corollary is useful for the correctness of our algorithm and it follows from Theorem \ref{UGB}, 

\begin{cor1}\label{neo}
Let $w$ be a primitive walk of a graph $G$. Then for the binomial
$B_w$ it holds $B_w\in \mathcal {U}_G$ if and only if the set of its
odd edges and the set of its even edges do not contain a cycle of
$w$.
\end{cor1}

\textbf{Proof.} Let $w$ be a primitive walk of a graph $G$. Since $w$ is primitive, from Corollary \ref{primitive-graph} it
follows that all of its blocks are either cycles either cut edges. From
Theorem~\ref{UGB}, $B_w\in \mathcal {U}_G$ if and only if $w$ is
mixed. By definition, $w$ is mixed if and only if all of its cyclic blocks
are not pure, i.e. for every block $B_i$ of its blocks,
$E(B_i)\nsubseteq E(\textbf{w}^+)$ and $E(B_i)\nsubseteq E(\textbf{w}^-)$, where by $E(B_i), E(\textbf{w}^+), E(\textbf{w}^-) $ we denote the sets of edges of the corresponding graphs. 

From Theorem \ref{UGB}, we are interested only in the mixed walks $w$ of the graph, which means that we are interested only in the cyclic blocks of the walk $w$. The result now follows. \hfill $\square$

Let see an example.

\begin{ex1}{\rm
We stay on Example \ref{exGraver} and we consider the walk $$w=\{e_1,e_2,e_3,e_4,e_5,e_6,e_7,e_8,e_9,e_{10},e_{11},e_{12}\},$$ where $B_w=\overbrace{e_1e_3e_5e_7e_9e_{11}}^{\textbf{w}^+}-\overbrace{e_2e_4e_6e_8e_{10}e_{12}}^{\textbf{w}^-}$ is its corresponding binomial. We note that the cycles of the walk $w$ are:
$$c_1=(e_1,e_2,e_3),c_2=(e_4,e_8,e_{12}),c_3=(e_5,e_6,e_7),c_4=(e_9,e_{10},e_{11})$$ We remark that there is the cycle
$c_2$ such that its edges are all even. By the previous corollary it follows that $B_w\notin \mathcal {U}_G$.}
\end{ex1}

Now we are ready to describe our algorithm. Its
correctness is guaranteed by Corollary \ref{neo}.
The algorithm takes as input the elements $B_w\in Gr_G$ as graphs $(V(B_w),E(B_w))$ separated into its two components, $\textbf{w}^+$ and $\textbf{w}^-$ in their graph representation (as, not-induced, subgraphs of $G$), $(V(\textbf{w}^+), E(\textbf{w}^+))$ and
$(V(\textbf{w}^-), E(\textbf{w}^-))$ respectively (i.e. the vertices and edges of each subgraph). We remark that the graph representation is not a multigraph, i.e. edges appear only once regardless of the number of times they appear in $B_w$. At each iteration, the algorithm considers one by one the elements $B_w\in Gr_G$ using their graph representation.
For each $B_w$, the goal of the algorithm is to check whether there is a cycle either in
$(V(\textbf{w}^+), E(\textbf{w}^+))$ either in $(V(\textbf{w}^-), E(\textbf{w}^-))$. According to Corollary~\ref{neo},
if such a cycle exists then $B_w \not \in \mathcal {U}_G$, otherwise $B_w \in \mathcal {U}_G$.

More specifically, by Corollary~\ref{primitive-graph}, the vertices of the graphs $(V(\textbf{w}^+), E(\textbf{w}^+))$ and
$(V(\textbf{w}^-), E(\textbf{w}^-))$ have degree either 1 either 2. The algorithm, thus, attempts to build a cycle first
in $(V(\textbf{w}^+), E(\textbf{w}^+))$  and then in $(V(\textbf{w}^-), E(\textbf{w}^-))$ by choosing as starting points vertices of degree 2, if such vertices exist. If such a cycle is found then, by Corollary~\ref{neo}, $B_w \not \in \mathcal {U}_G$, otherwise $B_w \in \mathcal {U}_G$.

The formal description of this algorithm, in pseudocode, follows below.

\begin{small}
\begin{tabbing}

\noindent {\bf ALGORITHM UGB} \\

\noindent {\bf INPUT}: \= \\
\> The elements $B_w\in Gr_G$, each split into its two parts, i.e.
$B_w^{+}$ and $B_w^{-}$, \\

\> given as simple (i.e. edges have multiplicity 1) subgraphs of $G$ \\
\> $(V(\textbf{w}^+), E(\textbf{w}^+))$  and $(V(\textbf{w}^-), E(\textbf{w}^-))$ respectively. \\

\noindent {\bf OUTPUT}: the elements of the
universal Gr\"{o}bner basis $\mathcal {U}_G$. \\

\noindent {\bf LOCAL} \= {\bf VARIABLES}:\\

\> Visited: \= $1 \times n$ array of boolean {\bf false} that is used in the cycle finding process. \\

\> stop: boolean variable. \\

\> $i, j$: integers. \\

\> $U$: a set of graph edges (i.e. non-ordered pairs of vertices). \\

1. {\bf begin algorithm} \\

2. $\mathcal {U}_G \leftarrow \emptyset$ \\

3. {\bf for} \= all elements $B_w\in Gr_G$ \\

4. \> stop $\leftarrow$ {\bf false} \\

5. \> {\bf for} $s \in \{+,-\}$ \\

6. \> \> Initialize {\em Visited} to {\bf false} \\

7. \> \> {\bf whi}\={\bf le} $\exists i \in V(\textbf{w}^s)$ such that $\deg(v_i) = 2$ {\bf and} Visited[$i$] = {\bf false} {\bf and} stop = {\em false}\\

8. \> \> \> Visited[$i$]  $\leftarrow$  {\bf true}. \\

9. \> \> \> $U \leftarrow \emptyset$ \\

10. \> \> \> {\bf rep}\={\bf eat} \\

11. \> \> \> \> Select $j$, such that $\{v_i, v_j\} \not \in U$ {\bf and} $\{v_i, v_j\} \in E(\textbf{w}^s)$\\

12. \> \> \> \> $U \leftarrow U \cup \{v_i, v_j\}$  \\

13. \> \> \> \> {\bf if} \= Visited[$j$] = {\bf true} \\

14. \> \> \> \> \> stop  $\leftarrow$  {\bf true} \\

15. \> \> \> {\bf until} \= stop = {\bf true} {\bf or} $\deg(v_j) = 1$ \\

16. \> \> \>  {\bf if} \= $\deg(v_j) = 1$ // No cycle was encountered during this iteration. \\

17. \> \> \> \> $E(\textbf{w}^s) \leftarrow E(\textbf{w}^s) - U$ // Delete encountered edges, updating vertices' degrees. \\

18. \> \> {\bf end while} // ``While'' at line 7. \\

19. \> {\bf end for} // ``For'' at line 5. \\

20.\>  {\bf if} \= stop = {\bf false} \\

21. \> \>  $\mathcal {U}_G \leftarrow \mathcal {U}_G \cup B_w$ \\

22. {\bf end for} // ``For'' at line 3. \\

23. {\bf return} $\mathcal {U}_G$ \\

24. {\bf end algorithm}

\end{tabbing}

\end{small}
\noindent
In the next theorem we prove the correctness and determine computational complexity of the
algorithm UGB and show that in order to decide algorithmically
whether a single element of the Graver basis of $I_G$ belongs also
to $\mathcal{U}_G$
it requires polynomial, in the number of the vertices of $G$, computational steps. By $|Gr_G| $
we denote the number of the elements of the Graver basis of $I_G$
and by $|B_w|$ the number of the edges in the walk $w$.
\par

\begin{thm1}\label{thm}
The Algorithm UGB correctly computes the set $\mathcal{U}_G$ for a
toric ideal of a given graph $G$ and its time
complexity (number of graph edges followed by the algorithm) is $$|Gr_G| \cdot \max_w{(O(|B_w|))}.$$
\end{thm1}

\noindent\textbf{Proof.} The algorithm is based on Corollary~\ref{neo}.
Accordingly, the loop at lines 3 to 22 iterates
over all elements $B_w$ in the Graver basis of the ideal.
For each such element $B_w$ the algorithm examines, in turn, $B_w^{+}$ and $B_w^{-}$ in
the loop at lines 5 to 19, for existence of cycles within either of these two sets.
The cycle finding process is in the loop at lines 7 to 18.
The ``stop'' variable is used in order to stop the iterations
when a cycle is found (it takes the value ``true'' and stops the iterations
at lines 7 to 18).

The loop at lines 7 to 18 selects, at each iteration,
an, yet, unvisited vertex of degree two, if it exists, and attempts to build a cycle starting from it. The cycle
building part is composed of the loop at lines 10 to 15.
The loop, repeatedly, selects consecutive edges until either a cycle
is found either a vertex of degree one is encountered. Note that according to Corollary~\ref{primitive-graph},
the vertices of the graphs that correspond to the elements $B_w$ of the Graver set have degree either one either two.

The ``if'' statement at lines 16 checks whether the loop at lines 10
to 15 was terminated due to finding a vertex of degree one, in which
case no cycle was detected. In this case, the encountered edges in
the set $U$ are discarded and the cycle detection process starts,
again, at line 7.

After the loop at lines 5 to 19 terminates, the ``if'' statement checks whether a cycle
was detected in either $B_w^{+}$ or $B_w^{-}$. This is the case where the variable
``stop'' has remained to its initial value {\bf false}. Thus, the element of the
Graver basis that was examined is inserted
into $\mathcal {U}_G$ otherwise the element is discarded. Then the algorithm starts over the iteration
at lines 3 to 22 to examine the next element of the Graver basis.
According to Corollary~\ref{neo} the algorithm correctly computes the elements of the Gr\"obner basis.

With respect to time complexity, there are $|Gr_G|$ iterations of the loop at lines 3 to 22.
For each such iteration, the graph formed by $B_w$ is examined at lines 5 to 19, separated into the components
$B_w^{+}$ and $B_w^{-}$ in subgraph representation.
Then the loop at lines 7 to 18 iteratively selects a vertex of degree 2 in order to follow consecutive
vertices, at the loop at lines 10 to 15, towards detecting a cycle, if it exists. Each vertex is considered once by the loop.
Then the loop at lines 10 to 15 follows consecutive edges in order
to detect a cycle or reach a vertex of degree 1. Again, the
encountered edges are considered once by the loop and once more, at
most, after the loop, at line 17.
Thus, in total, the complexity of the algorithm for a specific $B_w$ is in the order of $|V(B_w)| + |E(B_w)|$ or (in complexity theoretic
notation) $O(|V(B_w)| + |E(B_w)|)$. In our convention, the number of edges of the graph corresponding to $B_w$ is the number of its edges, denoted by $|B_w|$.
Thus, the complexity for $B_w$ is $O(|V(B_w)| + |B_w|)$. As the
number of vertices in a connected graph cannot be more than its edges minus 1,
the complexity is $O(|B_w|)$. Thus, in the worst case, the
complexity of the algorithm for locating all elements $B_w$ in the
Graver basis $Gr_G$ that also belong to the Gr\"obner basis is
$|Gr_G| \cdot  \max_w{(O(|B_w|))}$. \hfill $\square$
\
\

In~\cite{TT2} an upper bound for the
degrees of the primitive elements of $I_G$ is described. This bound is described
in the next proposition. It will be useful in estimating the
computational complexity of the algorithm.

\begin{prop1}\cite[Corollary
4.2.]{TT2} \label{b} Let $G$ be a graph with $n$ vertices, $n\geq 4$.
The largest degree $d$ of any binomial in the Graver basis (and in
the universal Gr\"{o}bner basis) for $I_{G}$ is $d\leq n-2$.
\end{prop1}

From the above proposition we know that $|B_w|\leq2(n-2)$. Therefore,
the following corollary holds, with $O(f(n))$ denoting functions of
$n$ bounded above by $cf(n)$, for some constant $c > 0$ and all $n >
n_0$, for some constant $n_0$.

\begin{cor1}
\label{cor-compl} The time complexity of Algorithm UGB is $O(|Gr_G|
\cdot n)$.
\end{cor1}

Let us consider two examples of the operation of the algorithm based
on Example~\ref{exGraver}. 

\begin{ex1}\rm{We first consider the operation of the
algorithm on the element $B_1\in Gr_G$ that appears in Figure~\ref{example} with
$B_1^{+}$ in blue color and $B_1^{-}$ in orange. The algorithm
enters the loop at line 7 trying to identify an unvisited vertex of
degree two. We assume that it selects $v_1$. Then it starts selecting
consecutive vertices at the loop at line 10. There is only one such
vertex, $v_3$, along edge $e_3$, which is followed by the algorithm.
Since the degree of $v_3$ is one, the algorithm stops and tries, at
line 7, to select another vertex of degree two, after discarding edge
$e_3$. The only choices are $v_4$ or $v_7$. The algorithm works as
with vertex $v_1$ concluding that in $B_1^{+}$ there is no cycle.
Then the loop at line 5 proceeds with $B_1^{-}$, i.e. the orange
edges. There are two choices of vertices of degree two. Let us assume
that at line 7 $v_1$ is chosen. Then it is clear that the loop at
line 10 will follow the vertices, detecting the cycle and rejecting
$B_1$.
}
\end{ex1}

\begin{ex1}
\rm{
Let us, now, turn to the example in Figure \ref{example1} with the
algorithm first considering the blue edges, i.e. $B_5^{+}$. 

\begin{figure}[h]
\begin{center}
\psfrag{A}{$e_1$}\psfrag{B}{$e_2$}\psfrag{C}{$e_3$}\psfrag{E}{$e_5$}\psfrag{F}{$e_6$}\psfrag{G}{$e_7$}\psfrag{H}{$e_8$}\psfrag{L}{$e_{12}$}
\psfrag{M}{$v_1$}\psfrag{N}{$v_2$}\psfrag{O}{$v_3$}\psfrag{P}{$v_4$}\psfrag{Q}{$v_5$}\psfrag{R}{$v_6$}\psfrag{S}{$v_7$}
\includegraphics{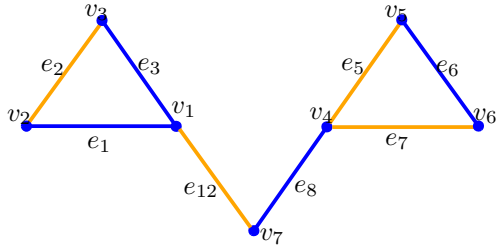}
\caption{The binomial $B_5=e_1e_3e_8^2e_6-e_2e_{12}^2e_7e_5$}
\label{example1}
\end{center}
\end{figure}

Again, the algorithm selects a vertex of degree two trying to locate a
cycle, e.g. vertex $v_1$. Then, the algorithm selects the next
unvisited vertex, say $v_3$. At this point, it has reached a vertex
of degree one and, thus, it stops without detecting a cycle. Then it
removes the edge it followed to repeat the process. However, now at
line 7 the algorithm does not find any vertex with degree two,
proceeding with $B_5^{-}$ (orange edges) at line 5. Now there is
only one vertex of degree two, $v_4$, which the algorithm selects. Let
us assume that it, then, selects $v_5$. This vertex has degree one
and, thus, the algorithm deletes the vertex and returns to line 7 in
order to select another vertex of degree two. However, such a vertex
does not exist and, thus, the algorithm is directed to line 21,
without having detected a cycle. Therefore, it includes the element
$B_5$ in $\mathcal {U}_G$.
}
\end{ex1}
\section{Conclusion}

We provide a new algorithm for computing the universal Gr\"obner basis of a toric ideal associated to incidence matrices of graphs.

As it is evident from the proof of Theorem~\ref{thm}, the algorithm
UGB decides in polynomial time whether a single element of the
Graver basis of $I_G$ belongs also to $\mathcal{U}_G$. 
However, since the algorithm must test exhaustively all elements of
the Graver basis for inclusion in the universal Gr\"{o}bner basis,
the factor $|Gr_G|$ is present in the final complexity figure. For
some classes of graphs this factor is not prohibitively large, i.e.
it is polynomial in $n$. For instance, the class of graphs that
consists of even cycles connected by paths where each cut vertex
belongs to exactly two blocks has at most $\frac{n}{4}$ cycles since
every cycle has at least four edges. Therefore its Graver basis has
cardinality $|Gr_G| \leq \frac{n}{4}$. For such graphs, the
complexity of the algorithm UGB is (from Corollary~\ref{cor-compl})
$O(n^2)$ which is a polynomial in the number $n$ of vertices of the
graph.\par
However, in many cases the size of the Graver basis appears to grow
exponentially fast with $n$. Thus, due to its large size, in
general, there can be no polynomial time algorithm for computing the
elements of the Graver basis of a general toric ideal. As an
indication of the computational difficulty of this problem,
in~\cite{LST} J.~De Loera, B.~Sturmfels and R.~Thomas showed,
computationally, that the number of the primitive elements of
$I_{K_8}$ is 45570, where $K_8$ is the complete graph on eight
vertices.\par
In general, there are no (to the best of our knowledge) general
tight upper bounds to the cardinality of the Graver basis of a graph
since this cardinality depends on its structure. Only a rough idea for the size of the Graver basis of a toric ideal
of a graph $G$ can be obtained from the bound to the degrees of its
elements (see Proposition~\ref{b}). However, the actual size of the
Graver basis itself can be huge.\par
For toric ideals of graphs, the problem of computing algorithmically
the Graver basis reduces to finding the set of even closed walks of
a corresponding graph $G$, which is again a very large set. We note that for the even closed walks which we are interested in, by Proposition \ref{b} we pass from each edge at most two times which means that the set of the even closed walks which we are looking for is finite. To the
best of our knowledge, existing algorithmic techniques can only
sample elements of the Graver basis in polynomial time
(see~\cite{OgTak}) or compute Graver basis by specialized
mathematical software, such as 4ti2 (see~\cite{ti}), for toric ideals of graphs with small Graver basis.
\\
\\

\end{document}